\documentclass[final,leqno,onethmnum,onefignum,onetabnum]{siamltex704}

\usepackage{amsmath}
\usepackage{amsfonts}
\usepackage[dvips]{epsfig,graphicx}
\usepackage{algorithmic}
\usepackage{algorithm}
\usepackage[usenames]{color}

\definecolor{DodgerBlue4}{rgb}{0.06,0.31,0.55}

\newcommand{\ignore}[1]{}

\newcommand{\comments}[1]{}
\date{}

\title{The cycle--convergence of restarted GMRES for normal matrices is sublinear}

\author{
E. Vecharynski
\thanks{
Department of Mathematical and Statistical Sciences,
University of Colorado Denver, Denver, CO 80217
(yaugen.vecharynski@ucdenver.edu, julien.langou@ucdenver.edu).}
\and
J. Langou\footnotemark[1]
}

\begin{document}

\setcounter{page}{1}
\maketitle

\begin{abstract}
We prove that the cycle--convergence of the restarted GMRES applied to a system
of linear equations with a normal coefficient matrix is sublinear.
\end{abstract}




\pagestyle{myheadings}
\thispagestyle{plain}
\markboth{E. VECHARYNSKI AND J. LANGOU}
{CYCLE-CONVERGENCE OF RESTARTED GMRES FOR NORMAL MATRICES}

\section{Introduction}
The \textit{generalized minimal residual method} (GMRES) was originally introduced by Saad and Schultz 
\cite{SaadSchultz:86} in 1986, and has become a popular method for
solving non-Hermitian systems of linear equations
\begin{equation}\label{sys}
 Ax=b,\quad A \in \mathbb{C}^{ n\times{n}},\quad b \in \mathbb{C}^n.
\end{equation}

GMRES is classified as a Krylov subspace (projection) iterative method.
At every new iteration $i$, GMRES constructs an approximation $x_i$ to the exact solution of 
(\ref{sys}), such that the $2$-norm of the corresponding residual vector $r_i = b - A x_i$ 
is minimized over the affine space $r_0 + A \mathcal{K}_i \left( A,r_0 \right)$, i.e. 
\begin{equation}\label{GMRES}
 r_i = \displaystyle \min_{u \in \mathcal{K}_i \left( A,r_0 \right)}\Vert r_0 - Au \Vert,
\end{equation}
where $\mathcal{K}_i \left( A,r_0 \right)$ is the $i$-dimensional Krylov subspace
\begin{equation*}
 \mathcal{K}_i \left( A,r_0 \right) = \mbox{span} \lbrace r_0, A r_0, \dots, A^{i-1} r_0 \rbrace
\end{equation*}
induced by the matrix $A$ and the initial residual vector $r_0 = b - A x_0$ with $x_0$ being
an initial approximate solution of (\ref{sys}).

As usual, in a linear setting, a notion of minimality is adjoint to some
orthogonality condition.  In our case, the minimization~(\ref{GMRES}) is
equivalent to forcing the new residual vector $r_i$ to be orthogonal to the
subspace $A\mathcal{K}_i \left( A,r_0 \right)$ (also known as the Krylov
residual subspace). In practice, for a large problem size, the latter
orthogonality condition results in a costly procedure of orthogonalization
against the expanding Krylov residual subspace.  Orthogonalization together
with storage requirement makes the GMRES method complexity and storage
prohibitive for practical application.  A straightforward treatment for this
complication is the so-called restarted GMRES~\cite{SaadSchultz:86}.

The \textit{restarted GMRES}, or GMRES($m$), is based on restarting GMRES after
every $m$ iterations.  At each restart, we use the latest approximate solution
as the initial approximation for the next GMRES run. Within this framework a
single run of $m$ GMRES iterations is called a GMRES($m$) cycle, $m$ is called
the restart parameter. Consequently, restarted GMRES can be regarded as a
sequence of GMRES($m$) cycles. When the convergence happens without any restart
occurring, the algorithm is known as \textit{full GMRES}.

Dealing with the restarted GMRES, our interest will shift towards the residual vectors $r_k$ at the
end of every $k$-th GMRES($m$) cycle (as opposed to the residual vectors $r_i$ (\ref{GMRES}) 
at each iteration of the original algorithm).
\begin{definition}[cycle--convergence] \label{def:1}
We define the cycle--convergence of restarted GMRES ($m$) as 
the norm of the residual vectors $\|r_k\|$ at the
end of every $k$-th GMRES($m$) cycle.
\end{definition}

We note that each $r_k$ satisfies the local minimality condition
\begin{equation}\label{GMRESm}
 r_k = \displaystyle \min_{u \in \mathcal{K}_m \left( A,r_{k-1} \right)}\Vert r_{k-1} - Au \Vert,
\end{equation}
where $\mathcal{K}_m \left( A,r_{k-1} \right)$ is the $m$-dimensional Krylov subspace produced at
the $k$-th GMRES($m$) cycle,
\begin{equation}\label{KrylovSubm}
 \mathcal{K}_m \left( A,r_{k-1} \right) = \mbox{span} \lbrace r_{k-1}, A r_{k-1}, \dots, A^{m-1} r_{k-1} \rbrace.
\end{equation}

The price paid for the reduction of the computational work, as follows from (\ref{GMRESm}) and (\ref{KrylovSubm}),
is the loss of global optimality (\ref{GMRES}). Although (\ref{GMRESm}) implies a monotonic decrease
of the norms of the residual vectors $r_k$, GMRES($m$) can stagnate \cite{SaadSchultz:86, ZavorinOLearyElman:03}.
This is in contrast with full GMRES which is guaranteed to converge to the exact 
solution of (\ref{sys}) in $n$ steps (assuming exact arithmetic). 
However, a proper choice of a preconditioner or/and a restart parameter, e.g.  
\cite{Embree:03, ErhelBurragePohl:96, Joubert:94}, can significantly accelerate 
the convergence of GMRES($m$), thus making the method practically attractive.

While a lot of efforts have been put into the characterization of the convergence of full
GMRES, e.g. \cite{VorstVuik:93, Eiermann:93, GreenbaumPtakStrakos:96, GreenbaumStrakos:94, Ipsen:00, Toh:97, 
LTrefethen:90}, 
our understanding of the behavior of GMRES($m$) is far from complete, leaving us with more questions than answers, 
e.g.~\cite{Embree:03}. In this manuscript, we prove that
the cycle--convergence of restarted 
GMRES for normal matrices
is sublinear.
This statement means that the reduction 
in the norm of the residual vector at the current GMRES($m$) cycle cannot be better than the
reduction at the previous cycle.

The current manuscript was inspired by ideas introduced in the technical report \cite{Zavorin:02} by I. Zavorin.  
In this work the author shows that, at every step of GMRES, a diagonalizable matrix $A$ and its Hermitian 
transpose $A^H$ yield the same worst-case behavior, and derives a necessary condition 
(the so-called cross-equality) for the worst-case right-hand side vector. We inherit the mathematical
tools for our analysis from \cite{Zavorin:02}, as well as \cite{Ipsen:00, ZavorinOLearyElman:03}, and
give their brief description, slightly adapted to the case of the restarted GMRES and a normal matrix $A$,
in Section 2. The main result of the sublinear cycle--convergence is proved in Section 3. In Section 4, the
behavior of GMRES($m$) in the nonnormal case is discussed.

\section{Krylov matrix, its pseudoinverse and spectral factorization}

Throughout the manuscript we will assume (unless otherwise explicitly stated) $A$ 
to be nonsingular and normal, i.e. $A$ allows the decomposition
\begin{equation}\label{nspectr}
 A=V \Lambda V^H,
\end{equation}
where $\Lambda \in \mathbb{C}^{n\times{n}}$ is a diagonal matrix with the diagonal elements
being the nonzero eigenvalues of $A$, and $V \in \mathbb{C}^{n\times{n}}$ is a unitary matrix 
of the corresponding eigenvectors.

Let us denote the $k$-th cycle of GMRES($m$) applied to the system (\ref{sys}) with the initial 
residual vector $r_{k-1}$ as GMRES($A$, $m$, $r_{k-1}$), $1 \leq m \leq n-1$. We assume that the 
residual vector $r_k$, produced at the end of GMRES($A$, $m$, $r_{k-1}$), is nonzero.

According to (\ref{GMRESm}) a run of GMRES($A$, $m$, $r_{k-1}$) entails the Krylov subspace 
$\mathcal{K}_m\left( A, r_{k-1} \right)$ (\ref{KrylovSubm}). For each  $\mathcal{K}_m\left( A, r_{k-1} \right)$
we define a matrix $K\left( A, r_{k-1} \right) \in \mathbb{C}^{n\times{(m+1)}}$, such that
\begin{equation}\label{KrylovMatr}
 K\left( A,r_{k-1} \right)=\left[r_{k-1} \quad Ar_{k-1} \quad  \dots \quad A^mr_{k-1}\right], \quad
k = 1, 2, \dots, q,
\end{equation}
where $q$ is the total number of GMRES($m$) cycles.

The matrix (\ref{KrylovMatr}) is called the Krylov matrix. We will say that
$K\left( A,r_{k-1} \right)$ corresponds to the cycle GMRES($A$, $m$, $r_{k-1}$). 
Note that the columns of $K\left( A, r_{k-1} \right)$ span the next, $(m+1)$-dimensional, 
Krylov subspace $\mathcal{K}_{m+1}(A,r_{k-1})$. By the assumption that $r_k \neq 0$,
\begin{equation*}
 \mbox{rank}\left(K\left( A, r_{k-1} \right)\right) = m + 1.
\end{equation*}

This latter equality allows us to introduce the Moore-Penrose pseudoinverse of the matrix $K\left( A, r_{k-1} \right)$,
\begin{equation*}
 K^{\dagger} \left( A, r_{k-1} \right) = 
\left( K^H \left( A,r_{k-1} \right) K\left( A,r_{k-1} \right) \right)^{-1} K^H \left( A,r_{k-1} \right)
\in \mathbb{C}^{{(m+1)} \times{n}},
\end{equation*}
which is well-defined and unique. The following lemma shows that the first column of 
$\left( K^{\dagger} \left( A, r_{k-1} \right) \right)^H$ is the next residual vector $r_k$
up to a scaling factor. 

\begin{lemma} \label{lem:1}
Given $A \in \mathbb{C}^{n\times{n}}$ (not necessarily normal) and the full rank Krylov matrix
$K\left( A, r_{k-1} \right) \in \mathbb{C}^{n\times{(m+1)}}$, corresponding to the cycle
GMRES($A$, $m$, $r_{k-1}$) for any $k = 1, 2, \dots, q$. Then
\begin{equation}\label{resk}
\left(K^{\dagger}\left( A,r_{k-1} \right) \right)^H e_1 = \frac{1}{\left\|r_k\right\|^2}r_k,  
\end{equation}
where $e_1 = \left[ 1 \quad 0 \quad \dots \quad 0 \right]^T \in \mathbb{R}^{m+1}$.
\end{lemma}

\textit{Proof}.
See Ipsen \cite[Theorem 2.1]{Ipsen:00}, as well as \cite{IpsenChandra:95, Stewart:87}.\hfill $\square$

Another important idea, mentioned in \cite{Ipsen:00} and intensively used in 
\cite{Zavorin:02, ZavorinOLearyElman:03}, provides the so-called spectral factorization
of the Krylov matrix $K\left( A,r_{k-1} \right)$ into three components, each one 
encapsulating separately the information on eigenvalues of $A$, its eigenvectors and 
the previous residual vector $r_{k-1}$. 

\begin{lemma} \label{lem:2}
Let $A \in \mathbb{C}^{n\times{n}}$ satisfying (\ref{nspectr}). Then the Krylov matrix 
$K\left( A,r_{k-1} \right)$, for any $k = 1, 2, \dots q$, can be factorized as  
\begin{equation}\label{KrylovSpectral}
 K\left( A,r_{k-1} \right) = V D_{k-1} Z,  
\end{equation}
where $d_{k-1} = V^H r_{k-1} \in \mathbb{C}^{n}$, $D_{k-1} = \mbox{diag}\left( d_{k-1} \right) \in \mathbb{C}^{n \times n}$ 
and $Z \in \mathbb{C}^{n \times {(m+1)}}$ is the Vandermonde matrix computed from the eigenvalues of $A$, 
\begin{equation}\label{vandermonde}
 Z = \left[ e \quad \Lambda e \quad \dots \quad \Lambda^m e \right],
\end{equation}
$e = \left[1 \quad 1 \quad \dots \quad 1 \right]^T \in \mathbb{R}^{n}$. 
\end{lemma}

\textit{Proof}.
Starting from (\ref{nspectr}) and the definition of the Krylov matrix (\ref{KrylovMatr})
\begin{eqnarray}
\nonumber K\left( A, r_{k-1} \right) & = & \left[ r_{k-1} \quad A r_{k-1} \quad \dots \quad A^m r_{k-1} \right]      ,\\
\nonumber   & = & \left[ V V^H r_{k-1} \quad V \Lambda V^H r_{k-1} \quad \dots \quad V \Lambda^m V^H r_{k-1} \right] ,\\
\nonumber   & = & V \left[ d_{k-1} \quad \Lambda d_{k-1} \quad \dots \quad \Lambda^m d_{k-1} \right]                 ,\\
\nonumber   & = & V \left[ D_{k-1} e \quad  \Lambda D_{k-1} e \quad \dots \quad  \Lambda^m D_{k-1} e \right]         ,\\
\nonumber   & = & V D_{k-1} \left[ e \quad  \Lambda e \quad \dots \quad \Lambda^m e \right] = V D_{k-1} Z.
\end{eqnarray}
\hfill$\square$

It is clear that the statement of Lemma~\ref{lem:2} can be easily generalized to the case of a diagonalizable (nonnormal)
matrix $A$ providing that we define $d_{k-1} = V^{-1} r_{k-1}$ in the lemma.

\section{The sublinear cycle--convergence of GMRES($m$)}

Along with (\ref{sys}) let us consider the system
\begin{equation}\label{sysH}
 A^H x = b
\end{equation}
with the matrix $A$ replaced by its Hermitian transpose. Clearly, according to (\ref{nspectr}),
\begin{equation}\label{nspectrH}
 A^H = V \overline{\Lambda} V^H.
\end{equation}

It turns out that $m$ steps of GMRES
applied to the systems (\ref{sys}) and (\ref{sysH}) produce the residual vectors of equal
norms, provided that the initial residual vector is the same for both GMRES runs. This observation
is crucial in concluding the sublinear cycle--convergence of GMRES($m$) and is formalized in the following
lemma.

\begin{lemma} \label{lem:3}
Let $r_m$ and $\hat{r}_m$ be the nonzero residual vectors obtained by applying $m$ steps of GMRES
to the systems (\ref{sys}) and (\ref{sysH}) respectively, $1 \leq m \leq n-1$. Then
\begin{equation*}
 \left\| r_m \right\| = \left\| \hat{r}_m \right\|,
\end{equation*}
provided that the initial approximate solutions of (\ref{sys}) and (\ref{sysH}) induce
the same initial residual vector $r_0$.
\end{lemma}

\textit{Proof}.
Consider a polynomial $p(z) \in \mathcal{P}_m$, where $\mathcal{P}_m$ is the set of all 
polynomials of degree at most $m$ defined on the complex plane, such that $p(0)=1$.
Let $r_0$ be a nonzero initial residual vector for the systems (\ref{sys}) and (\ref{sysH}) 
simultaneously. Since the matrix $A$ is normal, so is $p(A)$, thus $p(A)$ commutes with its
Hermitian transpose $p^H (A)$. We have 
\begin{equation*} 
\begin{split}
\| p(A) r_0 \|^2 & =  \langle p(A)r_0, p(A)r_0 \rangle
= \langle r_0, p^H(A) p(A) r_0 \rangle                                                                         ,\\
& = \langle r_0, p(A) p^H(A) r_0 \rangle 
= \langle p^H(A) r_0, p^H(A) r_0 \rangle                                                                       ,\\
& = \langle \left( V p(\Lambda) V^H \right)^H r_0, \left( V p(\Lambda) V^H \right)^H r_0 \rangle
= \langle V \overline{p}(\overline{\Lambda}) V^H r_0, V \overline{p}(\overline{\Lambda}) V^H r_0 \rangle       ,\\
& = \langle \overline{p}(V \overline{\Lambda} V^H) r_0, \overline{p}(V \overline{\Lambda} V^H) r_0 \rangle
= \| \overline{p}(V \overline{\Lambda} V^H) r_0 \|^2,
\end{split} 
\end{equation*}
where $\overline{p}(z)\in \mathcal{P}_m$ is the polynomial obtained from $p(z)$ by conjugating its coefficients.
By (\ref{nspectrH}) we conclude that
\begin{equation*}
 \| p(A) r_0 \| = \| \overline{p}(A^H) r_0 \|.
\end{equation*}

Since the last equality holds for any $p(z) \in \mathcal{P}_m$ it will also hold for the (GMRES) polynomial
$p_m(z)$, which minimizes $\| p(A) r_0 \|$ over $\mathcal{P}_m$. This polynomial exists and is unique
\cite[Theorem 2]{GreenbaumTrefthen:94}. 
Thus,
\begin{equation*}
\begin{split}
\| r_m \| & =  \min_{p \in \mathcal{P}_m} \| p(A) r_0 \| = \| p_m(A) r_0 \| = \| \overline{p}_m(A^H) r_0 \|    ,\\
& = \min_{p \in \mathcal{P}_m} \| \overline{p}(A^H) r_0 \| = \| \hat{r}_m \|,
\end{split} 
\end{equation*}
which proves the lemma. Moreover, we note that the two GMRES polynomials constructed after $m$ steps of GMRES
applied to (\ref{sys}) and (\ref{sysH}) with the same initial residual vector are the same up
to the complex conjugation of coefficients. 
\hfill$\square$

In the framework of the restarted GMRES Lemma~\ref{lem:3} suggests that the cycles GMRES($A$, $m$, $r_{k-1}$)
and GMRES($A^H$, $m$, $r_{k-1}$) result in the residual vectors $r_k$ and $\hat{r}_k$ of the
same norm.

So far we are ready to state the main theorem.

\begin{theorem}[The sublinear cycle--convergence of GMRES($m$)] \label{the:1}
Let $r_k$ be a sequence of nonzero residual vectors produced by GMRES($m$) applied
to the system (\ref{sys}) with a nonsingular normal matrix $A \in \mathbb{C}^{n \times n}$, 
$1 \leq m \leq n-1$. Then
\begin{equation}\label{sublin}
 \frac{\| r_{k} \|}{\| r_{k-1}  \|} \leq \frac{\| r_{k+1} \|}{\| r_k  \|}, \quad k = 1, \dots, q-1, 
\end{equation}
where $q$ is the total number of GMRES($m$) cycles.
\end{theorem}

\textit{Proof}.
Left multiplication of both parts of (\ref{resk}) by $K^H \left( A, r_{k-1} \right)$ leads to
\begin{equation*}
 e_1= \frac{1}{\| r_k \|^2}K^H \left( A,r_{k-1} \right) r_k.
\end{equation*}
By (\ref{KrylovSpectral}) in Lemma~\ref{lem:2}, we factorize the Krylov matrix $K \left( A,r_{k-1} \right)$ in the equality
above:
\begin{equation*}
\begin{split}
 e_1 & = \frac{1}{\| r_k \|^2} \left( V D_{k-1} Z \right)^H r_k 
= \frac{1}{\| r_k \|^2} Z^H \overline{D}_{k-1} V^H r_k,                                              \\
& = \frac{1}{\| r_k \|^2} Z^H \overline{D}_{k-1} d_k. 
\end{split}
\end{equation*}
%
Applying complex conjugation to this equality (and observing that $e_1$ is real), we get
\begin{equation*}
 e_1 = \frac{1}{\| r_k \|^2} Z^T D_{k-1} \overline{d}_k. 
\end{equation*}
According to the definition of $D_{k-1}$ in Lemma~\ref{lem:2}, $D_{k-1} \overline{d}_k = \overline{D}_k d_{k-1}$, thus
\begin{equation*}
 e_1 = \frac{1}{\| r_k \|^2} Z^T \overline{D}_k d_{k-1} 
= \frac{1}{\| r_k \|^2} \left( Z^T \overline{D}_k V^H \right) r_{k-1}.
\end{equation*}
From (\ref{KrylovSpectral}) and (\ref{nspectrH}) we notice that
\begin{equation*}
 Z^T \overline{D}_k V^H = \left( V D_k \overline{Z}  \right)^H = K^H\left( A^H, r_k \right), 
\end{equation*}
which leads to the following equality
\begin{equation}\label{undersys}
 e_1 = \frac{1}{\| r_k \|^2} K^H\left( A^H, r_k \right) r_{k-1}. 
\end{equation}
Considering the residual vector $r_{k-1}$ as a solution of the underdetermined system (\ref{undersys}),
we can represent the latter as 
\begin{equation}\label{solundersys}
 r_{k-1} = \| r_k \|^2 \left( K^H\left( A^H, r_k \right) \right)^{\dagger} e_1 + w_k, 
\end{equation}
where $w_k \in \mbox{null}\left( K^H\left( A^H, r_k \right)\right)$. Moreover, since 
\begin{equation*}
 w_k \perp \left( K^H\left( A^H, r_k \right) \right)^{\dagger} e_1, 
\end{equation*}
by the Pythagorean theorem we obtain
$$
\| r_{k-1} \|^2  = \| r_k \|^4 \|\left( K^H\left( A^H, r_k \right) \right)^{\dagger} e_1 \|^2 + \|w_k\|^2 ,
$$
 \text{now since} $\left( K^H\left( A^H, r_k \right) \right)^{\dagger} = \left( K^{\dagger}\left( A^H, r_k \right) \right)^H)$, \text{we get}  
\begin{align*}
\| r_{k-1} \|^2  & = \| r_k \|^4 \| \left( K^{\dagger}\left( A^H, r_k \right) \right)^H e_1 \|^2 + \|w_k\|^2,
&& \text{and then by (\ref{resk})}),                                                                                               \\
& = \frac{\| r_k \|^4}{\| \hat{r}_{k+1} \|^2 }  + \|w_k\|^2,   \\
& \geq  \frac{\| r_k \|^4}{\| \hat{r}_{k+1} \|^2 },
\end{align*}
where $\hat{r}_{k+1}$ is the residual vector at the end of the cycle GMRES($A^H$, $m$, $r_k$). Finally,
\begin{equation*}
\frac{\| r_k \|^2}{\| r_{k-1} \|^2}  
\leq \frac{\| r_k \|^2 \| \hat{r}_{k+1} \|^2 }{\| r_k \|^4} = 
\frac{\| \hat{r}_{k+1} \|^2 }{\| r_k \|^2}, 
\end{equation*}
so that 
\begin{equation}\label{resineq}
 \frac{\| r_k \|}{\| r_{k-1} \|} \leq \frac{\| \hat{r}_{k+1} \|}{\| r_k \|}. 
\end{equation}

By Lemma~\ref{lem:3}, the norm of the residual vector $\hat{r}_{k+1}$ at the end of the cycle GMRES($A^H$, $m$, $r_k$)
is equal to the norm of the residual vector $r_{k+1}$ at the end of the cycle GMRES($A$, $m$, $r_k$), which completes
the proof of the theorem.
\hfill$\square$

\begin{figure}
 \includegraphics[scale=0.37]{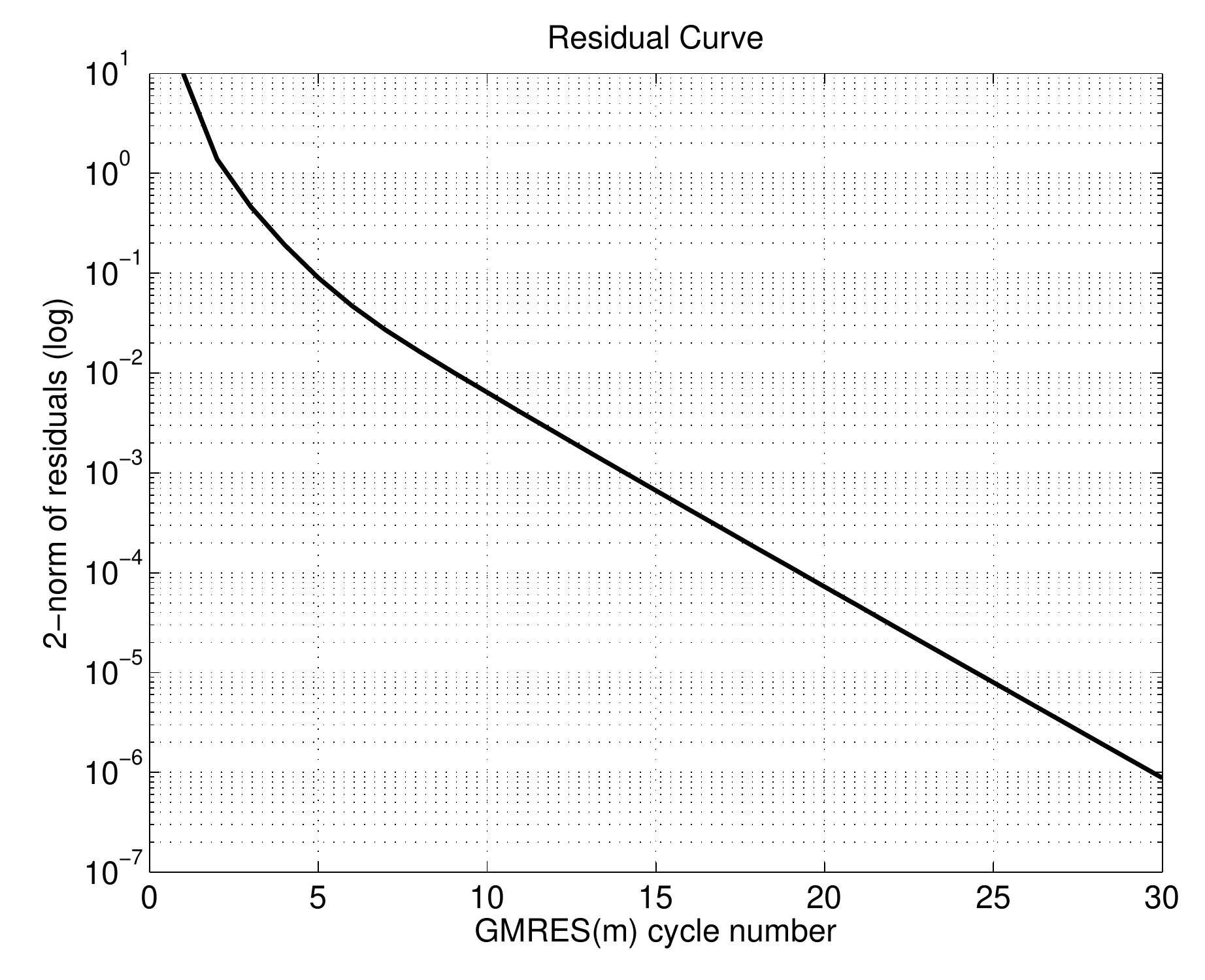}
 \includegraphics[scale=0.37]{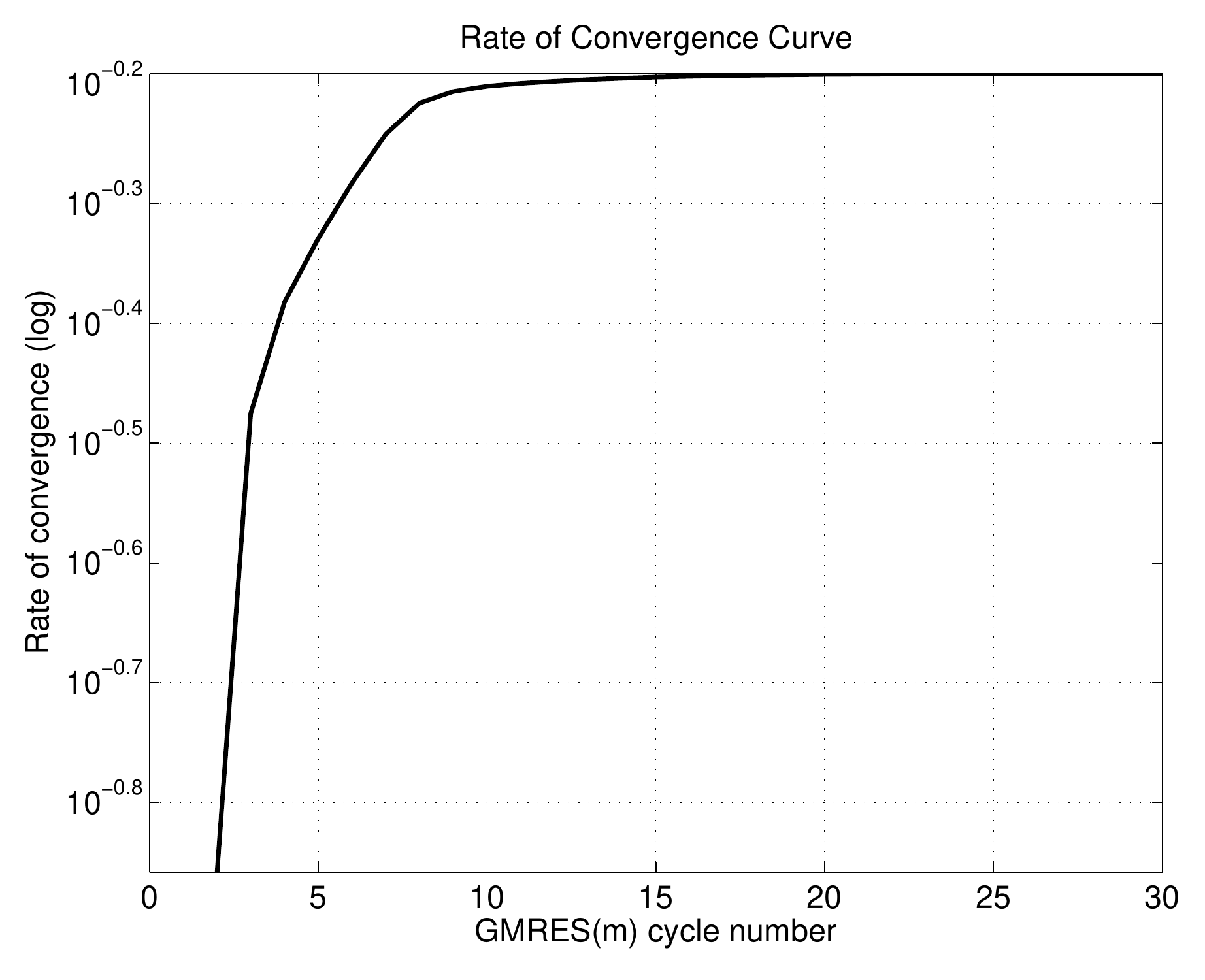}
 \caption{\label{fig:1}
 Cycle--convergence of GMRES(5) applied to a 100--by --100 normal matrix.
 }
\end{figure}

Geometrically, the theorem suggests that any residual curve of a restarted GMRES,
applied to a system with a nonsingular normal matrix, is nonincreasing and concave up (Figure~\ref{fig:1}).

From the proof of the Theorem~\ref{the:1} it is clear that, for a fixed $k$, the equality in (\ref{sublin}) holds if and only if
the vector $w_k$ (\ref{solundersys}) from the null space of the corresponding matrix $K^H\left( A^H, r_k \right)$ is
zero. In particular, when the restart parameter is chosen to be one less than the problem size, i.e. $m = n - 1$,
the matrix $K^H\left( A^H, r_k \right)$ in (\ref{undersys})
becomes
an $n$--by--$n$ nonsingular matrix, hence with a zero null space, and thus the Inequality~(\ref{sublin})
is indeed an equality when $m=n-1$.

It turns out that the cycle--convergence of GMRES($n-1$), applied to the system (\ref{sys}) 
with a nonsingular normal matrix $A$, can be completely determined by norms of the two initial
residual vectors $r_0$ and $r_1$.

\begin{corollary}[The cycle--convergence of GMRES($n-1$)] 
Given $\|r_0\|$ and $\|r_1\|$. Then, under assumptions of the Theorem~\ref{the:1}, norms of the residual vectors $r_k$ at the end
of each GMRES($n-1$) cycle obey the following formula
\begin{equation}\label{n-1}
\|r_{k+1}\| = \|r_1\| \left( \frac{\|r_1\|}{\|r_0\|}\right)^k, \quad k=1, \dots , q-1.
\end{equation}
\end{corollary}

\textit{Proof}.
The representation (\ref{solundersys}) of the residual vector $r_{k-1}$, for $m=n-1$, turns into
\begin{equation}\label{n-by-n_res}
 r_{k-1} = \| r_k \|^2 \left( K^H\left( A^H, r_k \right) \right)^{-1} e_1,
\end{equation}
implying, by the proof of the Theorem~\ref{the:1}, that the equality in (\ref{sublin}) holds at each GMRES($n-1$)
cycle. Thus,
\begin{equation*}
 \|r_{k+1}\| = \|r_k\| \frac{\|r_k\|}{\|r_{k-1}\|}, \quad k=1, \dots , q-1. 
\end{equation*}

We show (\ref{n-1}) by induction in $k$. Using the formula above, it is easy to verify (\ref{n-1}) for
$\| r_2 \|$ and $\| r_3 \|$ ($k = 1, 2$). Let's assume that for some $k$, $3 \leq k \leq q-1$, 
$\| r_{k-1} \|$ and $\| r_k \|$ can also be computed by (\ref{n-1}). Then
\begin{align*}
 \|r_{k+1}\| & = \|r_k\| \frac{\|r_k\|}{\|r_{k-1}\|} 
= \|r_1\| \left( \frac{\|r_1\|}{\|r_0\|}\right)^{k-1} 
\frac{\|r_1\| \left( \frac{\|r_1\|}{\|r_0\|}\right)^{k-1}}{\|r_1\| \left( \frac{\|r_1\|}{\|r_0\|}\right)^{k-2}} \\
& = \|r_1\| \left( \frac{\|r_1\|}{\|r_0\|}\right)^{k-1} \left( \frac{\|r_1\|}{\|r_0\|}\right)
= \|r_1\| \left( \frac{\|r_1\|}{\|r_0\|}\right)^k.
\end{align*}
Thus, (\ref{n-1}) holds for all $k=1, \dots , q-1$. 
\hfill$\square$

Another observation in the proof of the Theorem~\ref{the:1} leads to a well known result 
due to Baker, Jessup and Manteuffel~\cite{BakerJessupManteuffel:05}.
In this paper, the authors prove that,
when GMRES($n-1$) is applied to a system
with Hermitian or skew-Hermitian matrix,
the residual vectors at the end of each restart cycle alternate direction in a cyclic fashion~\cite[Theorem 2]{BakerJessupManteuffel:05}.
In the following corollary we (slightly) refine this result by providing the exact expression 
for the constants $\alpha_k$ in~\cite[Theorem 2]{BakerJessupManteuffel:05}.  

\begin{corollary}[The alternating residuals]
Let $r_k$ be a sequence of nonzero residual vectors produced by GMRES($n-1$) applied
to the system (\ref{sys}) with a nonsingular Hermitian or skew-Hermitian matrix $A \in \mathbb{C}^{n \times n}$. 
Then
\begin{equation}\label{alt}
r_{k+1}=\alpha_k r_{k-1}, \quad \alpha_k = \frac{\left\|r_{k+1}\right\|^2}{\left\|r_k\right\|^2} \in \left(0,1\right], 
\quad k=1,2, \dots, q-1. 
\end{equation}
\end{corollary}

\textit{Proof}.
For the case of a Hermitian matrix $A$, i.e. $A^H = A$, the proof follows directly from 
(\ref{n-by-n_res}) and (\ref{resk}).

Let $A$ be skew-Hermitian, i.e. $A^H = -A$. Then, by (\ref{n-by-n_res}) and (\ref{resk}),
\begin{equation*}
 r_{k-1} =  \left( K^H\left( A^H, r_k \right) \right)^{-1} e_1
=  \left( K^H\left( -A, r_k \right) \right)^{-1} e_1
= \frac{\| r_k \|^2}{\| \hat{r}_{k+1} \|^2} \hat{r}_{k+1},
\end{equation*}
where $\hat{r}_{k+1}$ is the residual vector produced at the end of the cycle GMRES($-A$, $n-1$, $r_k$).

According to (\ref{GMRESm}), the residual vectors $r_{k+1}$ and $\hat{r}_{k+1}$ at the end of the cycles
GMRES($A$, $n-1$, $r_k$) and GMRES($-A$, $n-1$, $r_k$) are obtained by orthogonalizing $r_k$ against 
the Krylov residual subspaces $A\mathcal{K}_{n-1} \left( A, r_k \right)$ and 
$\left(-A\right)\mathcal{K}_{n-1} \left( -A, r_k \right)$ respectively.
But $\left(-A\right)\mathcal{K}_{n-1} \left( -A, r_k \right) = A\mathcal{K}_{n-1} \left( A, r_k \right)$, hence
$\hat{r}_{k+1} = r_{k+1}$.
\hfill$\square$

\section{Note on the departure from normality}

\begin{figure}
 \includegraphics[scale=0.37]{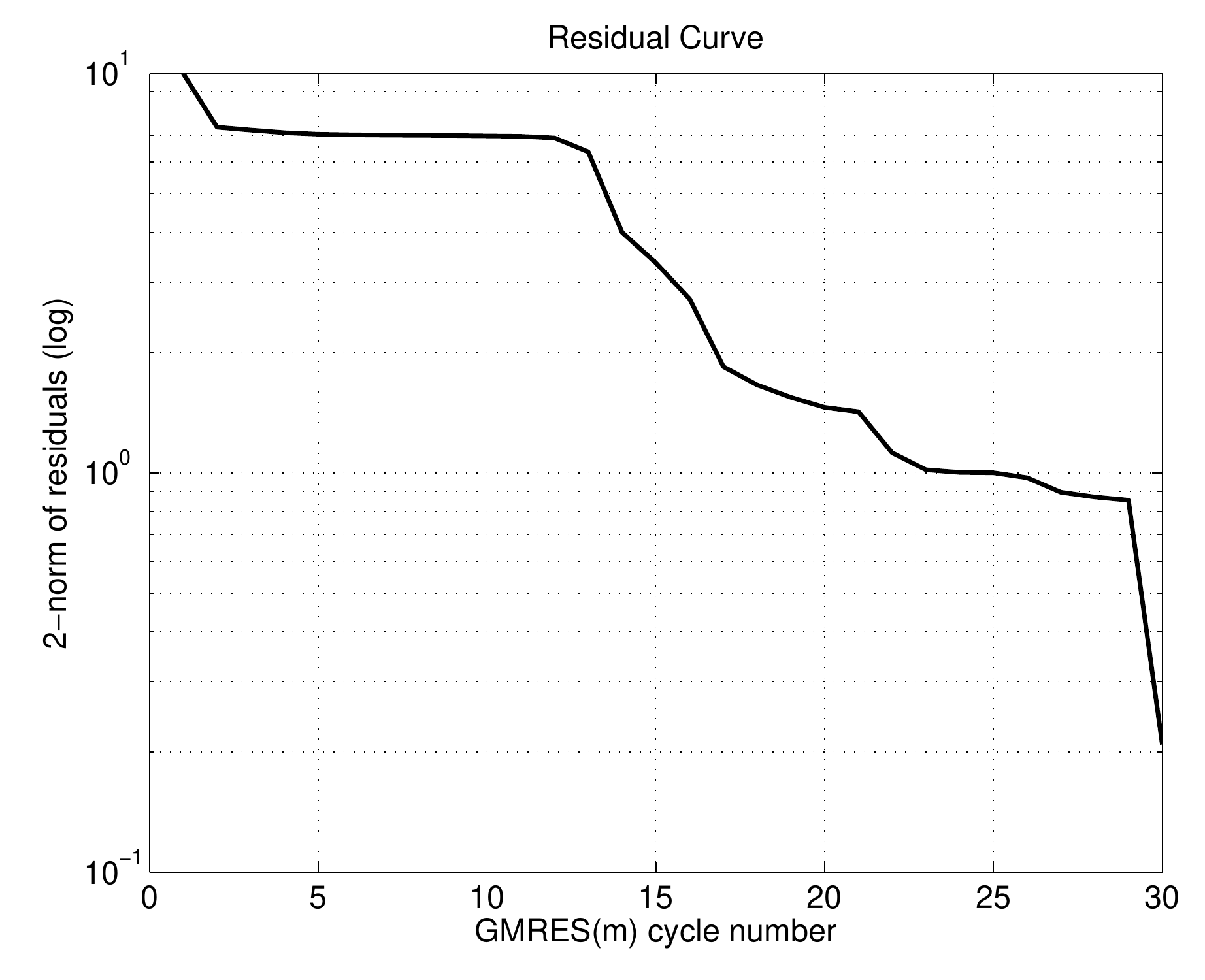}
 \includegraphics[scale=0.37]{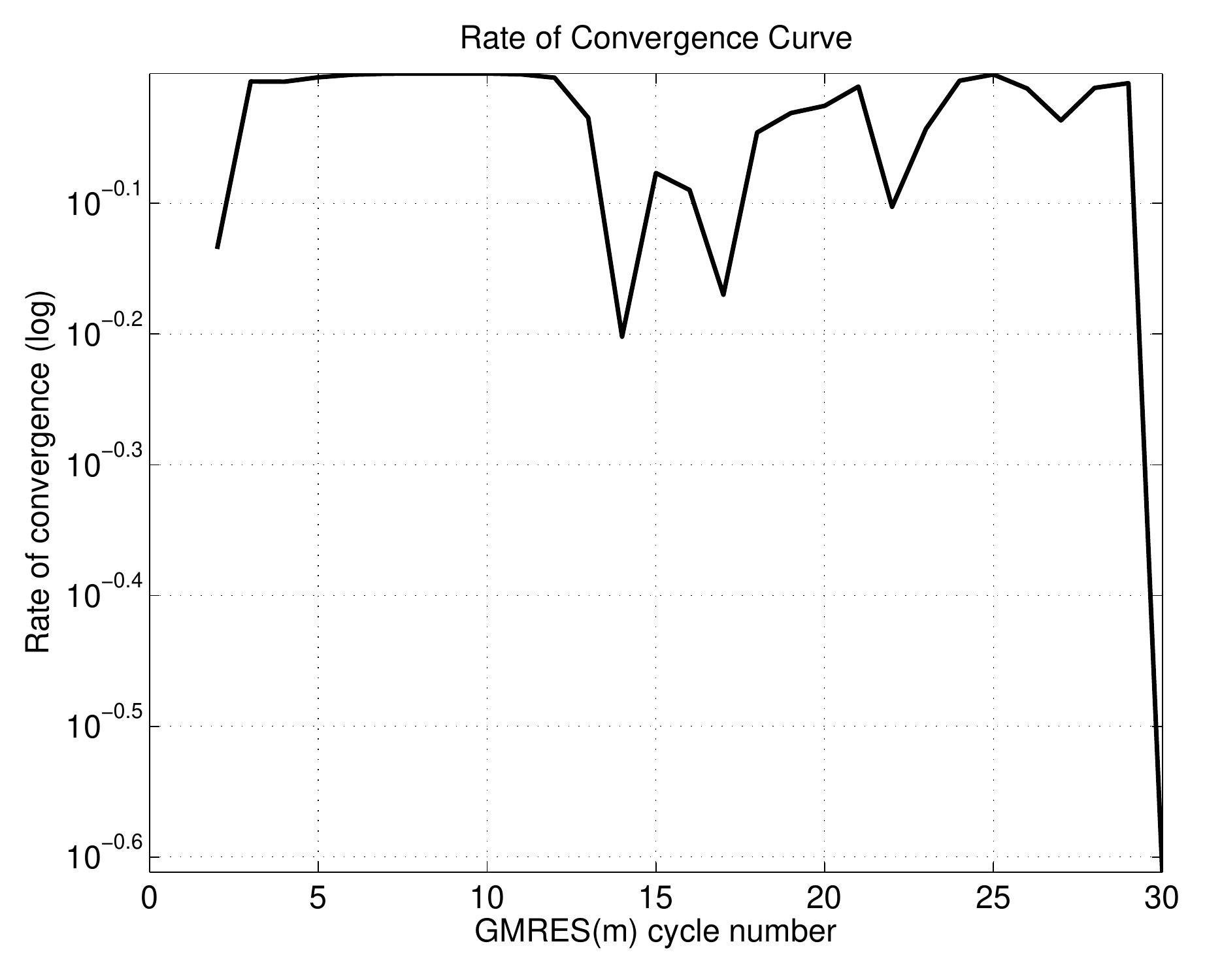}
 \caption{\label{fig:2}
 Cycle--convergence of GMRES(5) applied to a 100--by--100 diagonalizable (nonnormal) matrix.
 }
\end{figure}

In general, for systems with nonnormal matrices, the cycle--convergence
behavior of the restarted GMRES is not sublinear. In Figure~\ref{fig:2}, we
take a nonnormal diagonalizable matrix for illustration purpose and one can
observe the claim. Indeed, for
nonnormal matrices, it has been
observed the cycle--convergence of restarted GMRES can be
superlinear~\cite{ZhongMorgan:08}.

In this concluding section we restrict our attention to the case of a diagonalizable matrix $A$,
\begin{equation}\label{spectrd}
 A = V \Lambda V^{-1}, \quad A^H = V^{-H} \overline{\Lambda} V^H.
\end{equation}

The analysis performed in Theorem~\ref{the:1} can be generalized
for the case of a diagonalizable matrix (\cite{Zavorin:02}), resulting in the inequality (\ref{resineq}).
However, as we depart from normality, Lemma~\ref{lem:3} fails to hold and the norm of the residual vector
$\hat{r}_{k+1}$ at the end of the cycle GMRES($A^H$, $m$, $r_k$) is no longer equal to the norm of the
vector $r_{k+1}$ at the end of GMRES($A$, $m$, $r_k$). Moreover, since the eigenvectors of $A$ can be
significantly changed by the Hermitian conjugation, as (\ref{spectrd}) suggests, the matrices
$A$ and $A^H$ can have almost nothing in common, so that the norms of $\hat{r}_{k+1}$ and $r_{k+1}$ are,
possibly, far from being equal. This
gives a chance for breaking
the sublinear convergence of GMRES($m$), provided that the subspace $A\mathcal{K}_m \left( A,r_k \right)$
results in a better approximation (\ref{GMRESm}) of the vector $r_k$ than the subspace
$A^H \mathcal{K}_m \left( A^H ,r_k \right)$.

It is natural to expect that the convergence of the restarted GMRES for ``almost normal'' matrices will be ``almost sublinear''.
We quantify this statement in the following lemma.

\begin{lemma} \label{lem:4}
Let $r_k$ be a sequence of nonzero residual vectors produced by GMRES($m$) applied
to the system (\ref{sys}) with a nonsingular diagonalizable (\ref{spectrd}) matrix 
$A \in \mathbb{C}^{n \times n}$, $1 \leq m \leq n-1$. Then
\begin{equation}\label{asublin}
\frac{\| r_{k} \|}{\| r_{k-1} \|} \leq \frac{\alpha \left( \| r_{k+1} \| + \beta_k \right)}{\| r_k \|}, \quad k = 1, \dots, q-1, 
\end{equation}
where $\alpha = \frac{1}{\sigma_{min}^2(V)}$, $\beta_k = \|p_k(A) (I-VV^H)r_k\|$,
$p_k(z)$ is the polynomial constructed at the cycle GMRES($A$, $m$, $r_k$), and
where $q$ is the total number of GMRES($m$) cycles.
Note that as $V^HV\longrightarrow I$,
$0 < \alpha\longrightarrow 1$ and $0 < \beta_k \longrightarrow 0$.
\end{lemma}

\textit{Proof}.
Consider the norm of the residual vector $\hat{r}_{k+1}$ at the end of the cycle GMRES($A^H$, $m$, $r_k$).
\begin{equation*}
 \| \hat r_{k+1} \| = \min_{\hat{p} \in \mathcal{P}_m} \| \hat{p}(A^H) r_k \| \leq \| p(A^H) r_k \|, 
\end{equation*}
where $p(z) \in \mathcal{P}_m$ is any polynomial of degree at most $m$, such that $p(0)=1$. Then,
using (\ref{spectrd}), 
\begin{eqnarray}
\nonumber \| \hat r_{k+1} \| & \leq & \| p(A^H) r_k \|\\
\nonumber                    &  =   & \| V^{-H} p(\overline{\Lambda}) V^H r_k \|\\
\nonumber                    &  =   & \| V^{-H} p(\overline{\Lambda})(V^{-1}V) V^H r_k \|\\
\nonumber                    &  =   & \| V^{-H} p(\overline{\Lambda}) V^{-1} (VV^H) r_k \|\\
\nonumber                    &  =   & \| V^{-H} p(\overline{\Lambda}) V^{-1} (I-(I-VV^H)) r_k \|\\
\nonumber                    &  =   & \| V^{-H} p(\overline{\Lambda})\left(V^{-1}r_k-V^{-1}(I-VV^H)r_k \right)\|\\
\nonumber                    & \leq & \| V^{-H} \| \|p(\overline{\Lambda})\left(V^{-1}r_k-V^{-1}(I-VV^H)r_k\right)\|.
\end{eqnarray}
Note that 
\begin{equation*}
\|p(\overline{\Lambda})\left(V^{-1}r_k-V^{-1}(I-VV^H)r_k\right)\|  =   \|\overline{p}(\Lambda)\left(V^{-1}r_k-V^{-1}(I-VV^H)r_k\right)\|.
\end{equation*}
Thus,
\begin{eqnarray}
\nonumber \| \hat r_{k+1} \| & \leq & \| V^{-H} \| \|\overline{p}(\Lambda)\left(V^{-1}r_k-V^{-1}(I-VV^H)r_k\right)\|\\
\nonumber                    &  =   & \| V^{-H} \| \|(V^{-1}V)\overline{p}(\Lambda)\left(V^{-1}r_k-V^{-1}(I-VV^H)r_k\right)\|\\
\nonumber                    & \leq & \| V^{-H} \|\|V^{-1}\| \|V\overline{p}(\Lambda)V^{-1}r_k - V\overline{p}(\Lambda)V^{-1}(I-VV^H)r_k\|\\
\nonumber                    &  =   & \frac{1}{\sigma_{min}^2(V)} \|\overline{p}(V\Lambda V^{-1})r_k - \overline{p}(V \Lambda V^{-1}) (I-VV^H)r_k\|\\
\nonumber                    & \leq & \frac{1}{\sigma_{min}^2(V)} \left(\|\overline{p}(A)r_k\| + \| \overline{p}(A) (I-VV^H)r_k\| \right),
\end{eqnarray}
where $\sigma_{min}$ is the smallest singular values of $V$.

Since the last inequality holds for any polynomial $\overline{p}(z) \in \mathcal{P}_m$, it will also hold for
$\overline{p}(z) = p_k(z)$, where $p_k(z)$ is the polynomial constructed at the cycle GMRES($A$, $m$, $r_k$). Hence,
\begin{equation*}
\nonumber \| \hat r_{k+1} \| \leq \frac{1}{\sigma_{min}^2(V)} \left(\|r_{k+1}\| + \|p_k(A) (I-VV^H)r_k\| \right).
\end{equation*}
Setting $\alpha = \frac{1}{\sigma_{min}^2(V)}$, $\beta_k = \|p_k(A) (I-VV^H)r_k\|$ and observing that
$\alpha\longrightarrow 1$, $\beta_k \longrightarrow 0$ as $V^HV\longrightarrow I$, from (\ref{resineq}),
we obtain (\ref{asublin}).
\hfill$\square$


\bibliography{GMRES_sublinear}

\end{document}